\magnification 1200
\input plainenc
\input amssym
\fontencoding{T2A}
\inputencoding{utf-8}
\tolerance 4000
\relpenalty 10000
\binoppenalty 10000
\parindent 1.5em

\hsize 17truecm
\vsize 23.5truecm
\hoffset 0truecm
\voffset -0.5truecm

\font\TITLE labx1440
\font\tenrm larm1000
\font\cmtenrm cmr10
\font\tenit lati1000
\font\tenbf labx1000
\font\tentt latt1000
\font\teni cmmi10 \skewchar\teni '177
\font\tensy cmsy10 \skewchar\tensy '60
\font\tenex cmex10
\font\teneufm eufm10
\font\eightrm larm0800
\font\cmeightrm cmr8
\font\eightit lati0800
\font\eightbf labx0800
\font\eighttt latt0800
\font\eighti cmmi8 \skewchar\eighti '177
\font\eightsy cmsy8 \skewchar\eightsy '60
\font\eightex cmex8
\font\eighteufm eufm8

\font\cmsixrm cmr6

\font\sixbf labx0600
\font\sixi cmmi6 \skewchar\sixi '177
\font\sixsy cmsy6 \skewchar\sixsy '60
\font\sixeufm eufm6

\font\cmfiverm cmr5

\font\fivebf labx0500
\font\fivei cmmi5 \skewchar\fivei '177
\font\fivesy cmsy5 \skewchar\fivesy '60
\font\fiveeufm eufm5
\font\tencmmib cmmib10 \skewchar\tencmmib '177
\font\eightcmmib cmmib8 \skewchar\eightcmmib '177
\font\sevencmmib cmmib7 \skewchar\sevencmmib '177
\font\sixcmmib cmmib6 \skewchar\sixcmmib '177
\font\fivecmmib cmmib5 \skewchar\fivecmmib '177
\newfam\cmmibfam
\textfont\cmmibfam\tencmmib \scriptfont\cmmibfam\sevencmmib
\scriptscriptfont\cmmibfam\fivecmmib
\def\tenpoint{\def\rm{\fam0\tenrm}\def\it{\fam\itfam\tenit}%
	\def\bf{\fam\bffam\tenbf}\def\tt{\fam\ttfam\tentt}%
	\textfont0\cmtenrm \scriptfont0\cmsevenrm \scriptscriptfont0\cmfiverm
  	\textfont1\teni \scriptfont1\seveni \scriptscriptfont1\fivei
  	\textfont2\tensy \scriptfont2\sevensy \scriptscriptfont2\fivesy
  	\textfont3\tenex \scriptfont3\tenex \scriptscriptfont3\tenex
  	\textfont\itfam\tenit
	\textfont\bffam\tenbf \scriptfont\bffam\sevenbf
	\scriptscriptfont\bffam\fivebf
	\textfont\eufmfam\teneufm \scriptfont\eufmfam\seveneufm
	\scriptscriptfont\eufmfam\fiveeufm
	\textfont\cmmibfam\tencmmib \scriptfont\cmmibfam\sevencmmib
	\scriptscriptfont\cmmibfam\fivecmmib
	\normalbaselineskip 12pt
	\setbox\strutbox\hbox{\vrule height8.5pt depth3.5pt width0pt}%
	\normalbaselines\rm}
\def\eightpoint{\def\rm{\fam 0\eightrm}\def\it{\fam\itfam\eightit}%
	\def\bf{\fam\bffam\eightbf}\def\tt{\fam\ttfam\eighttt}%
	\textfont0\cmeightrm \scriptfont0\cmsixrm \scriptscriptfont0\cmfiverm
	\textfont1\eighti \scriptfont1\sixi \scriptscriptfont1\fivei
	\textfont2\eightsy \scriptfont2\sixsy \scriptscriptfont2\fivesy
	\textfont3\eightex \scriptfont3\eightex \scriptscriptfont3\eightex
	\textfont\itfam\eightit
	\textfont\bffam\eightbf \scriptfont\bffam\sixbf
	\scriptscriptfont\bffam\fivebf
	\textfont\eufmfam\eighteufm \scriptfont\eufmfam\sixeufm
	\scriptscriptfont\eufmfam\fiveeufm
	\textfont\cmmibfam\eightcmmib \scriptfont\cmmibfam\sixcmmib
	\scriptscriptfont\cmmibfam\fivecmmib
	\normalbaselineskip 11pt
	\abovedisplayskip 5pt
	\belowdisplayskip 5pt
	\setbox\strutbox\hbox{\vrule height7pt depth2pt width0pt}%
	\normalbaselines\rm
}

\def\No{\char 157}
\def\empty{}

\catcode`\@ 11
\catcode`\" 13
\def"#1{\ifx#1<\char 190\relax\else\ifx#1>\char 191\relax\else #1\fi\fi}

\def\newl@bel#1#2{\expandafter\def\csname l@#1\endcsname{#2}}
\openin 11\jobname .aux
\ifeof 11
	\closein 11\relax
\else
	\closein 11
	\input \jobname .aux
	\relax
\fi

\newcount\c@section
\newcount\c@subsection
\newcount\c@subsubsection
\newcount\c@equation
\newcount\c@bibl
\c@section 0
\c@subsection 0
\c@subsubsection 0
\c@equation 0
\c@bibl 0
\def\lab@l{}
\def\label#1{\immediate\write 11{\string\newl@bel{#1}{\lab@l}}%
	\ifhmode\unskip\fi}

\def\section#1{\global\advance\c@section 1
	{\par\vskip 3ex plus 0.5ex minus 0.1ex
	\rightskip 0pt plus 1fill\leftskip 0pt plus 1fill\noindent
	{\bf\S\thinspace\number\c@section .~#1}\par\penalty 25000%
	\vskip 1ex plus 0.25ex}
	\gdef\lab@l{\number\c@section.}
	\c@subsection 0
	\c@subsubsection 0
	\c@equation 0
}
\def\subsection{\global\advance\c@subsection 1
	\par\vskip 1ex plus 0.1ex minus 0.05ex{\bf\number\c@subsection. }%
	\gdef\lab@l{\number\c@section.\number\c@subsection}%
	\c@subsubsection 0\c@equation 0%
}
\def\subsubsection{\global\advance\c@subsubsection 1
	\par\vskip 1ex plus 0.1ex minus 0.05ex%
	{\bf\number\c@subsection.\number\c@subsubsection. }%
	\gdef\lab@l{\number\c@section.\number\c@subsection.%
		\number\c@subsubsection}%
}
\def\equation{\global\advance\c@equation 1
	\gdef\lab@l{\number\c@section.\number\c@subsection.%
	\number\c@equation}{\rm\number\c@equation}
}
\def\bibitem#1{\global\advance\c@bibl 1
	[\number\c@bibl]%
	\gdef\lab@l{\number\c@bibl}\label{#1}
}
\def\ref@ref#1.#2:{\def\REF@{#2}\ifx\REF@\empty{\S\thinspace#1}%
	\else\ifnum #1=\c@section {#2}\else {\S\thinspace#1.#2}\fi\fi
}
\def\ref@eqref#1.#2.#3:{\ifnum #1=\c@section\ifnum #2=\c@subsection
	{(#3)}\else{#2\thinspace(#3)}\fi\else{\S\thinspace#1.#2\thinspace(#3)}\fi
}
\def\ref#1{\expandafter\ifx\csname l@#1\endcsname\relax
	{\bf ??}\else\edef\mur@{\csname l@#1\endcsname :}%
	{\expandafter\ref@ref\mur@}\fi
}
\def\eqref#1{\expandafter\ifx\csname l@#1\endcsname\relax
	{(\bf ??)}\else\edef\mur@{\csname l@#1\endcsname :}%
	{\expandafter\ref@eqref\mur@}\fi
}
\def\cite#1{\expandafter\ifx\csname l@#1\endcsname\relax
	{\bf ??}\else\hbox{\bf\csname l@#1\endcsname}\fi
}

\catcode`\@ 12

\immediate\openout 11\jobname.aux


\frenchspacing\rm
\leftline{УДК~517.927}\vskip 0.25truecm
{\leftskip 0cm plus 1fill\rightskip 0cm plus 1fill\parindent 0cm\baselineskip 15pt
\TITLE Об одной априорной мажоранте собственных значений задач Штурма--Лиувилля\par
\vskip 0.25truecm\rm \ А.$\,$А.~Владимиров
\par\bigskip
$$
	\abovedisplayskip 0pt
	\vbox{\hsize 0.75\hsize\leftskip 0cm\rightskip 0cm
	\eightpoint\rm
	{\bf Аннотация:\/} Изучается вопрос о точной априорной мажоранте $M_\gamma$
	наимень\-шего соб\-ственного значения задачи Штурма--Лиувилля
	$$
		\openup -1\jot
		\displaylines{-y''+qy=\lambda y,\cr y(0)=y(1)=0}
	$$
	с ограничениями на потенциал вида $q\leqslant 0$ и $\int_0^1 |q|^\gamma\,dx=1$,
	где $\gamma\in (0,1/2)$. Показывается, что эта мажоранта подчиняется строгой
	оценке $M_\gamma<\pi^2$.
	}
$$
}

\vskip 0.25truecm
\subsection
Рассмотрим граничную задачу
$$
	\displaylines{-y''+qy=\lambda y,\cr y(0)=y(1)=0,}
$$
где потенциал выбирается внутри семейства
$$
	A_\gamma=\left\{q\in C[0,1]\;:\;q\leqslant 0,\quad\int_0^1 |q|^\gamma\,dx=1\right\}.
$$
Свяжем с указанной задачей априорную мажоранту $M_\gamma\rightleftharpoons\sup_{q\in A_\gamma}
\lambda_0(q)$ наименьшего собственного значения. Как было установлено в [\cite{OSZE}:
Теорема~1.2], при $\gamma\geqslant 1/2$ выполняется равенство $M_\gamma=\pi^2$,
а при $\gamma<1/3$ справедлива строгая оценка $M_\gamma<\pi^2$. Данный результат
повторён также в ряде позднейших публикаций того же автора (см., например, [\cite{OZMF}:
Теорема~2.1]). Целью настоящей заметки является установление справедливости строгой оценки
$M_\gamma<\pi^2$ также в случае $\gamma\in [1/3,1/2)$.

\subsection
Далее мы всегда предполагаем величину $\gamma\in (0,1/2)$ каким-либо образом
зафиксированной. Пусть $q\in C[0,1]$~--- неположительный потенциал, удовлетворяющий
соотношению $\lambda_0(q)>(\pi-\varepsilon)^2$, где $\varepsilon>0$~--- некоторое
достаточно малое число. Рассмотрим связанные с соответствующей собственной функцией
$y\in C^2[0,1]$ функции $\varrho,\vartheta\in C^1[0,1]$ и $\sigma\in C[0,1]$ вида
$$
	\varrho\cdot\pmatrix{\sin\vartheta\cr\cos\vartheta}
		\rightleftharpoons\pmatrix{y\cr y'/\sqrt{\lambda_0(q)}},\qquad
		\sigma\rightleftharpoons |q|\,\sin^2\vartheta.
$$
Функция $\vartheta$ подчиняется уравнению
$$
	\vartheta'={1\over\varrho}\cdot\pmatrix{\cos\vartheta&
		-\sin\vartheta}\cdot\pmatrix{y'\cr y''/\sqrt{\lambda_0(q)}}
		={\lambda_0(q)+\sigma\over\sqrt{\lambda_0(q)}}
$$
и пробегает, строго возрастая, отрезок $[0,\pi]$. Соотношения
$$
	\int_0^1{\sigma\vartheta'\over\lambda_0(q)+\sigma}\,dx=
		\int_0^1\left[\vartheta'-\sqrt{\lambda_0(q)}\right]\,dx
		=\pi-\sqrt{\lambda_0(q)}<\varepsilon
$$
означают, что множество
$$
	E_\varepsilon\rightleftharpoons\left\{x\in [0,1]\;:\;
		\sigma(x)>\varepsilon^{(1-2\gamma)/(1-\gamma)}\right\}
$$
подчиняется оценке
$$
	\int_{E_\varepsilon}\vartheta'\,dx<\mu(\varepsilon)\rightleftharpoons
		\pi^2\varepsilon^{\gamma/(1-\gamma)}+\varepsilon.
$$
Далее мы всегда будем предполагать выполненным неравенство $\mu(\varepsilon)<\pi$.
В этом случае справедливы оценки
$$
	\displaylines{\int_{E_\varepsilon}\sin^{-2\gamma}\vartheta\cdot
		\vartheta'\,dx\leqslant 2\int_0^{\mu(\varepsilon)/2}\sin^{-2\gamma}x\,dx
		<{2\pi^2\over 1-2\gamma}\,\varepsilon^{(1-2\gamma)\gamma/(1-\gamma)},\cr
		\int_{\overline{E_\varepsilon}}\sin^{-2\gamma}\vartheta\cdot
		\vartheta'\,dx\leqslant\int_0^\pi\sin^{-2\gamma}x\,dx<{4\over 1-2\gamma}.
	}
$$
С учётом заведомо выполняющихся, ввиду $\lambda_0(q)>4$, неравенств
$$
	2\sqrt{\lambda_0(q)}\sigma^\gamma(x)\leqslant\lambda_0(q)+\sigma(x)
$$
это означает справедливость соотношений
$$
	\eqalign{\int_0^1 |q|^{\gamma}\,dx&=\sqrt{\lambda_0(q)}\cdot\int_0^1
		{\sigma^{\gamma}\sin^{-2\gamma}\vartheta\over\lambda_0(q)+\sigma}\cdot
		\vartheta'\,dx\cr &=\sqrt{\lambda_0(q)}\cdot\left[
		\int_{E_\varepsilon}{\sigma^{\gamma}\sin^{-2\gamma}
		\vartheta\over\lambda_0(q)+\sigma}\cdot\vartheta'\,dx+
		\int_{\overline{E_\varepsilon}}{\sigma^{\gamma}\sin^{-2\gamma}
		\vartheta\over\lambda_0(q)+\sigma}\cdot\vartheta'\,dx\right]\cr
		&\leqslant{1\over 2}\cdot\int_{E_\varepsilon}\sin^{-2\gamma}
		\vartheta\cdot\vartheta'\,dx+{\varepsilon^{(1-2\gamma)\gamma/(1-\gamma)}
		\over\sqrt{\lambda_0(q)}}\cdot\int_{\overline{E_\varepsilon}}
		\sin^{-2\gamma}\vartheta\cdot\vartheta'\,dx\cr
		&\leqslant{\pi^2+2\over	1-2\gamma}\cdot
		\varepsilon^{(1-2\gamma)\gamma/(1-\gamma)}.
	}
$$
Тем самым, справедливость соотношения $q\in A_\gamma$ несовместима с предположением
о значительной малости параметра $\varepsilon>0$. Это и означает выполнение искомой
оценки $M_{\gamma}<\pi^2$.

\vskip 0.4cm
\eightpoint\rm
{\leftskip 0cm\rightskip 0cm plus 1fill\parindent 0cm
\bf Литература\par\penalty 20000}\vskip 0.4cm\penalty 20000
\bibitem{OSZE} {\it С.$\,$С.~Ежак\/}. Оценки первого собственного значения задачи
Штурма--Лиувилля с условиями Дирихле~/ В~кн.: Качественные свойства решений
дифференциальных уравнений и смежные вопросы спектрального анализа. М.:~ЮНИТИ-ДАНА,
2012. С.~517--559.

\bibitem{OZMF} {\it С.$\,$С.~Ежак\/}. Об одной задаче минимизации функционала,
порождённого задачей Штурма--Лиувилля с интегральным условием на потенциал~//
Вестник СамГУ.~--- 2015.~--- \No~6$\,$(128).~--- С.~57--61.

\bye